\documentclass[11pt]{amsart}
\usepackage[margin=1in]{geometry}
\usepackage{amsmath,amssymb,amsthm,mathtools}
\usepackage{enumitem}
\usepackage{booktabs}
\usepackage[colorlinks=true,linkcolor=blue,citecolor=blue]{hyperref}

\newtheorem{theorem}{Theorem}[section]
\newtheorem{lemma}[theorem]{Lemma}
\newtheorem{proposition}[theorem]{Proposition}
\newtheorem{corollary}[theorem]{Corollary}
\newtheorem{mainthm}{Theorem}

\theoremstyle{definition}
\newtheorem{definition}[theorem]{Definition}
\newtheorem{example}[theorem]{Example}
\newtheorem{problem}[theorem]{Problem}
\theoremstyle{remark}
\newtheorem{remark}[theorem]{Remark}

\newcommand{\R}{\mathbb{R}}
\newcommand{\N}{\mathbb{N}}
\newcommand{\Z}{\mathbb{Z}}
\newcommand{\eps}{\varepsilon}
\newcommand{\ru}{\mathrm{ru}}
\newcommand{\othree}{\mathrm{o}_3}
\newcommand{\ocv}{\xrightarrow{\,\othree\,}}
\newcommand{\id}{\circ}
\newcommand{\CDf}{\textup{(CD$^\flat$)}}
\newcommand{\CD}{\textup{(CD)}}
\newcommand{\CDs}{\textup{(CD$_\sigma$)}}

\begin{document}

\title[Pre-Riesz-valued vector metrics I]{Vector metrics with values in pre-Riesz spaces I:\\ transfer of order convergence along vector lattice covers}

\author{C\"uneyt \c{C}ev\.{\i}k}
\address{Department of Mathematics, Faculty of Science, Gazi University, Ankara, Turkey}
\email{ccevik@gazi.edu.tr}


\begin{abstract}
Vector metric spaces in the sense of \c{C}evik and Altun take values in a vector lattice, but their axioms use only the order of the value space. We develop the theory for vector metrics with values in an Archimedean pre-Riesz space $E$, i.e.\ an ordered vector space admitting a bipositive order dense embedding $i$ into a vector lattice $Y$. We ask when an $E$-valued vector metric $d$ and the vector-lattice-valued metric $i\id d$ have the same convergent and Cauchy nets and sequences.

Relatively uniform convergence transfers along every cover. For nets, $\othree$-convergence transfers if and only if the cover has an intrinsic domination property. This property holds for pervasive covers, by a result of Kalauch and van Gaans, for uniformly dense covers and in finite dimensions; it fails for an example of Kalauch and van Gaans. Sequential convergence transfers under a countable domination property implied by uniform density, but not under pervasiveness. For uncountable $\Gamma$, the vector lattice $\R\mathbf 1\oplus c_{00}(\Gamma)$ is pervasive in its Dedekind completion $\ell^\infty(\Gamma)$ without the transfer, while for countable $\Gamma$ the transfer holds; the difference is explained by order separability. In finite dimensions, with a closed generating cone, all these convergences coincide with norm convergence in every cover.

For $C^1[0,1]$-valued integral metrics, positivity of the kernel on the diagonal forces all order convergences to coincide with $L^1$-convergence. For the Volterra kernel $(t-s)_+$, by contrast, $L^1$-convergence, relatively uniform convergence and order convergence are pairwise different, and we describe each explicitly.
\end{abstract}

\maketitle

\section{Introduction}\label{sec:intro}

Vector metric spaces were introduced by \c{C}evik and Altun~\cite{CA09}: the distance takes its values in a vector lattice $E$, and a sequence converges if its distances to the limit are dominated by a sequence decreasing to $0$ in $E$ (Definition~\ref{a:def:vm}). Continuity and fixed point theorems for such spaces were studied in \cite{AC11,Cev14}. An older relative is Kantorovich's theory of $K$-metric spaces \cite{Zab97}, where convergence is governed by a norm or a lattice structure on the value space.

The axioms of a vector metric use only the order of $E$. A lattice structure is needed only once moduli, suprema and infima enter the convergence theory. Pre-Riesz spaces \cite{vH93,KvG19} are exactly the ordered vector spaces that embed bipositively and order densely into a vector lattice. Every directed Archimedean ordered vector space is pre-Riesz; typical examples are $C^1[0,1]$, spaces of polynomials, and $\R^n$ with a closed generating non-simplicial cone. For an $E$-valued vector metric $d$ and a vector lattice cover $(Y,i)$, the map $i\id d$ is a classical vector metric, and the question is when convergence transfers back from $Y$ to $E$. Contractions and fixed points are studied in the companion paper \cite{VC2}, and completions in \cite{DC}.

We now state the main results. Let $E$ be an Archimedean pre-Riesz space with vector lattice cover $(Y,i)$, and $d$ an $E$-valued vector metric.

\begin{mainthm}[Theorem~\ref{d:thm:charact}, Theorem~\ref{a:thm:findim}]\label{thmB}
For nets, $\othree$-convergence transfers from $i\id d$ to $d$ for all $E$-valued vector metrics if and only if the cover has the domination property \CDf\ of Definition~\ref{d:def:CD}. If $\dim E<\infty$ and the cone is closed and generating, then all order convergences coincide with norm convergence of $d(x_\lambda,x)$ in every cover, which need be neither pervasive nor uniformly dense.
\end{mainthm}

\begin{mainthm}[Proposition~\ref{a:prop:ru}, Theorems~\ref{a:thm:o3} and \ref{a:thm:seq}]\label{thmA}
Relatively uniform convergence, Cauchy nets and completeness of $d$ and $i\id d$ coincide for every cover. The $\othree$-convergences of $d$ and $i\id d$ coincide if $i(E)$ is pervasive or uniformly dense in $Y$, and sequential convergence in the sense of \cite{CA09} coincides if $i(E)$ is uniformly dense.
\end{mainthm}

The hypotheses cannot be dropped, and nets and sequences behave differently.

\begin{mainthm}[Theorem~\ref{a:thm:counter}, Example~\ref{d:ex:counter}]\label{thmC}
Without hypotheses the transfer fails, even for sequences. Pervasiveness does not suffice for sequences: for uncountable $\Gamma$, the space $\R\mathbf 1\oplus c_{00}(\Gamma)$ is a pervasive sublattice of its Dedekind completion $\ell^\infty(\Gamma)$ in which a sequence of unit vectors converges for $i\id d$ but not for $d$. For countable $\Gamma$ the transfer holds.
\end{mainthm}

Finally, for integral metrics with values in $C^1[0,1]$ the convergences can be compared explicitly.

\begin{mainthm}[Theorems~\ref{d:thm:band} and \ref{d:thm:volterra}]\label{thmD}
For integral metrics $d_k(f,g)(t)=\int_0^1k(t,s)|f(s)-g(s)|\,ds$ with values in $C^1[0,1]$, positivity of $k$ on the diagonal forces $d_k$-$\othree$, $d_k$-$\ru$ and $L^1$-convergence to coincide, even if $k$ vanishes off a band. For the Volterra kernel $(t-s)_+$ the three are pairwise different: $\ru$-convergence is convergence in $L^1((1-s)\,ds)$, and sequential order convergence is local $L^1$-convergence on $[0,1)$ together with boundedness in $L^1((1-s)\,ds)$.
\end{mainthm}

Theorem~\ref{thmB} and the uniformly dense, sequential and integral-metric results in Theorems~\ref{thmA}, \ref{thmC} and \ref{thmD} are new. For pervasive covers and the order convergence of $E$ itself, the transfer of nets is due to Kalauch and van Gaans \cite[Proposition~3.7.12]{KvG19}, and the space behind the failure without hypotheses in Theorem~\ref{thmC} is their example \cite[Example~3.7.10]{KvG19}.

The paper is organised as follows. Sections~\ref{a:sec:prelim} and \ref{a:sec:vm} contain preliminaries on pre-Riesz spaces, order convergences and vector metrics. Section~\ref{a:sec:transfer} proves the transfer results under (P) and (D), Section~\ref{a:sec:findim} treats finite-dimensional value spaces, and Section~\ref{a:sec:examples} gives examples and the counterexample. Section~\ref{d:sec:nets} characterises the transfer for nets, Section~\ref{d:sec:seq} treats sequences, Section~\ref{d:sec:integral} integral metrics, and Section~\ref{sec:open} lists open problems.

\section{Preliminaries}\label{a:sec:prelim}

Throughout, $E$ is a real ordered vector space with cone $E_+$, so that $E_++E_+\subseteq E_+$, $\alpha E_+\subseteq E_+$ for $\alpha\ge 0$, and $E_+\cap(-E_+)=\{0\}$. We write $x\le y$ if $y-x\in E_+$. For $A\subseteq E$ we write $A^u$ and $A^l$ for the sets of upper and lower bounds of $A$. All nets are indexed by directed sets, that is, nonempty preordered sets in which any two elements have a common upper bound (antisymmetry is not required), and a property holds \emph{eventually} for a net $(x_\lambda)_{\lambda\in\Lambda}$ if there is $\lambda_0$ such that it holds for all $\lambda\ge\lambda_0$. For a net $(a_\gamma)$ we write $a_\gamma\downarrow 0$ if it is decreasing and $\inf_\gamma a_\gamma=0$.

\begin{definition}\label{a:def:basic}
The space $E$ is \emph{Archimedean} if $nx\le y$ for all $n\in\N$ implies $x\le 0$. It is \emph{directed} if $E=E_+-E_+$. It is a \emph{pre-Riesz space} if for all $x,y,z\in E$ the inclusion $\{x+z,y+z\}^u\subseteq\{x,y\}^u$ implies $z\in E_+$.
\end{definition}

Every directed Archimedean ordered vector space is pre-Riesz, and every pre-Riesz space is directed \cite{vH93}, \cite[Section~2.4]{KvG19}. An element $u\in E_+$ is an \emph{order unit} of $E$ if for every $x\in E$ there is $\lambda\ge 0$ with $\pm x\le\lambda u$; in a vector lattice an order unit is also called a \emph{strong unit}. A subspace $Z$ of a vector lattice $Y$ is \emph{majorising} if every $y\in Y$ satisfies $y\le z$ for some $z\in Z$.

\begin{definition}\label{a:def:cover}
Let $Y$ be a vector lattice and $i\colon E\to Y$ linear.
\begin{enumerate}[label=(\alph*)]
\item The map $i$ is \emph{bipositive} if $i(x)\ge 0\Leftrightarrow x\ge 0$.
\item The subspace $i(E)$ is \emph{order dense} in $Y$ if $y=\inf\{i(x):x\in E,\ y\le i(x)\}$ for every $y\in Y$.
\item The pair $(Y,i)$ is a \emph{vector lattice cover} of $E$ if $i$ is bipositive and $i(E)$ is order dense in $Y$.
\end{enumerate}
\end{definition}

By van Haandel's theorem \cite{vH93}, \cite[Section~2.4]{KvG19}, $E$ has a vector lattice cover if and only if $E$ is pre-Riesz. Applying (b) to $-y$ shows that order density is equivalent to $y=\sup\{i(x):x\in E,\ i(x)\le y\}$ for all $y\in Y$. A bipositive linear map is injective.

\begin{lemma}\label{a:lem:major}
Let $(Y,i)$ be a vector lattice cover of $E$. For every $y\in Y$ there is $x\in E$ with $y\le i(x)$.
\end{lemma}

\begin{proof}
The set $\{i(x):x\in E,\ y\le i(x)\}$ has infimum $y$, so it is nonempty.
\end{proof}

Infima of subsets of $E$ can be computed in a cover.

\begin{lemma}\label{a:lem:inf}
Let $(Y,i)$ be a vector lattice cover of $E$, $A\subseteq E$ and $a\in E$. Then $\inf_E A=a$ if and only if $\inf_Y i(A)=i(a)$. In particular, a net $(a_\gamma)$ in $E$ satisfies $a_\gamma\downarrow 0$ in $E$ if and only if $i(a_\gamma)\downarrow 0$ in $Y$.
\end{lemma}

\begin{proof}
Let $\inf_E A=a$. Then $i(a)$ is a lower bound of $i(A)$. If $y\in Y$ is another lower bound, then every $x\in E$ with $i(x)\le y$ is a lower bound of $A$ by bipositivity, so $x\le a$. By order density, $y=\sup\{i(x):i(x)\le y\}\le i(a)$. Conversely, if $\inf_Y i(A)=i(a)$ and $z\in E$ is a lower bound of $A$, then $i(z)\le i(a)$ and hence $z\le a$.
\end{proof}

The Archimedean property passes to covers.

\begin{lemma}\label{a:lem:arch}
If $E$ is Archimedean and $(Y,i)$ is a vector lattice cover of $E$, then $Y$ is Archimedean.
\end{lemma}

\begin{proof}
Let $y,z\in Y$ with $ny\le z$ for all $n$. By Lemma~\ref{a:lem:major} choose $w\in E$ with $z\le i(w)$. If $x\in E$ and $i(x)\le y$, then $i(nx)\le z\le i(w)$, so $nx\le w$ for all $n$, and $x\le 0$. Hence $y=\sup\{i(x):i(x)\le y\}\le 0$.
\end{proof}

We shall use two conditions on a cover.

\begin{definition}\label{a:def:PD}
Let $(Y,i)$ be a vector lattice cover of $E$.
\begin{enumerate}[label=(\alph*)]
\item[(P)] The subspace $i(E)$ is \emph{pervasive} in $Y$ if for every $y\in Y$ with $y>0$ there is $x\in E$ with $0<i(x)\le y$ \cite[Section~2.8]{KvG19}.
\item[(D)] The subspace $i(E)$ is \emph{uniformly dense} in $Y$ if there is $u\in E_+$ such that for every $y\in Y$ and every $\eps>0$ there is $x\in E$ with $y\le i(x)\le y+\eps\, i(u)$.
\end{enumerate}
\end{definition}

If $i(u)$ is a strong unit of $Y$, then (D) says exactly that $i(E)$ is dense in $Y$ for the norm $\|\cdot\|_{i(u)}$: given $x$ with $\|i(x)-y\|_{i(u)}\le\eps/2$, the element $x+\tfrac{\eps}{2}u$ satisfies $y\le i(x+\frac\eps2 u)\le y+\eps\, i(u)$. Examples in Section~\ref{a:sec:examples} show that (P) and (D) are independent.

\begin{definition}\label{a:def:conv}
Let $(x_\lambda)_{\lambda\in\Lambda}$ be a net in $E$ and $x\in E$.
\begin{enumerate}[label=(\alph*)]
\item $x_\lambda\xrightarrow{\mathrm o}x$ if there is a net $(a_\lambda)_{\lambda\in\Lambda}$ in $E$ with $a_\lambda\downarrow 0$ and $\pm(x_\lambda-x)\le a_\lambda$ for all $\lambda$.
\item $x_\lambda\xrightarrow{\othree}x$ if there is a net $(a_\gamma)_{\gamma\in\Gamma}$ in $E$ with $a_\gamma\downarrow 0$ such that for every $\gamma$ eventually $\pm(x_\lambda-x)\le a_\gamma$.
\item $x_\lambda\xrightarrow{\ru}x$ if there is $u\in E_+$ such that for every $\eps>0$ eventually $\pm(x_\lambda-x)\le\eps u$.
\end{enumerate}
\end{definition}

These notions go back to Abramovich and Sirotkin \cite{AS05} and Wulich \cite{Wul17}; see also \cite[Section~1.7]{KvG19}. The essential difference between (a) and (b) is that in (b) the dominating net may have its own index set.

\section{Pre-Riesz-valued vector metrics}\label{a:sec:vm}

The axioms of a vector metric involve only the order of the value space.

\begin{definition}\label{a:def:vm}
Let $X$ be a nonempty set and $E$ a pre-Riesz space. A map $d\colon X\times X\to E$ is an \emph{$E$-valued vector metric} if for all $x,y,z\in X$
\begin{enumerate}[label=\textup{(VM\arabic*)}]
\item $d(x,y)=0$ if and only if $x=y$;
\item $d(x,y)\le d(x,z)+d(y,z)$.
\end{enumerate}
We call $(X,d,E)$ a \emph{pre-Riesz-valued vector metric space}.
\end{definition}

When $E$ is a vector lattice, this is the definition of \cite{CA09}.

\begin{lemma}\label{a:lem:vmbasic}
For all $x,y,z,x',y'\in X$:
\begin{enumerate}[label=\textup{(\roman*)}]
\item $d(x,y)\ge 0$ and $d(x,y)=d(y,x)$;
\item $\pm\bigl(d(x,z)-d(y,z)\bigr)\le d(x,y)$;
\item $\pm\bigl(d(x,y)-d(x',y')\bigr)\le d(x,x')+d(y,y')$.
\end{enumerate}
\end{lemma}

\begin{proof}
(i) Taking $y=x$ in (VM2) gives $0\le 2d(x,z)$. Taking $z=x$ gives $d(x,y)\le d(y,x)$, and interchanging $x$ and $y$ gives equality. (ii) By (VM2) and (i), $d(x,z)\le d(x,y)+d(z,y)=d(x,y)+d(y,z)$; interchange $x$ and $y$. (iii) By (ii), $d(x,y)\le d(x,x')+d(x',y)\le d(x,x')+d(x',y')+d(y',y)$, and symmetrically.
\end{proof}

In a vector lattice, (ii) reads $|d(x,z)-d(y,z)|\le d(x,y)$. The point of the $\pm$ form is that no modulus is needed.

\begin{definition}\label{a:def:vmconv}
Let $(x_\lambda)_{\lambda\in\Lambda}$ be a net in $X$ and $x\in X$.
\begin{enumerate}[label=(\alph*)]
\item $x_\lambda\to x$ in \emph{$d$-o} if $d(x_\lambda,x)\le a_\lambda$ for a net $(a_\lambda)_{\lambda\in\Lambda}$ in $E$ with $a_\lambda\downarrow 0$.
\item $x_\lambda\to x$ in \emph{$d$-$\othree$} if there is a net $(a_\gamma)_{\gamma\in\Gamma}$ in $E$ with $a_\gamma\downarrow 0$ such that for every $\gamma$ eventually $d(x_\lambda,x)\le a_\gamma$.
\item $x_\lambda\to x$ in \emph{$d$-$\ru$} if there is $u\in E_+$ such that for every $\eps>0$ eventually $d(x_\lambda,x)\le\eps u$.
\item A sequence $(x_n)$ converges to $x$ in the \emph{sequential sense} if $d(x_n,x)\le a_n$ for a sequence $a_n\downarrow 0$ in $E$ \cite{CA09}.
\end{enumerate}
A net is \emph{$d$-$\othree$-Cauchy} if there is $(a_\gamma)$ in $E$ with $a_\gamma\downarrow 0$ such that for every $\gamma$ eventually $d(x_\lambda,x_\mu)\le a_\gamma$ (that is, for all $\lambda,\mu\ge\lambda_0(\gamma)$). It is \emph{$d$-$\ru$-Cauchy} if there is $u\in E_+$ such that for every $\eps>0$ eventually $d(x_\lambda,x_\mu)\le\eps u$. A sequence is \emph{Cauchy in the sequential sense} if $d(x_n,x_{n+p})\le a_n$ for all $n,p$ and some sequence $a_n\downarrow 0$ in $E$. For each of these three kinds, $X$ is \emph{complete} if every Cauchy net (respectively sequence) converges in the same sense.
\end{definition}

The three kinds of convergence are related as follows.

\begin{proposition}\label{a:prop:relations}
Let $(X,d,E)$ be a pre-Riesz-valued vector metric space.
\begin{enumerate}[label=\textup{(\roman*)}]
\item Sequential convergence implies $d$-o convergence, and $d$-o convergence implies $d$-$\othree$ convergence. If $E$ is Archimedean, then $d$-$\ru$ convergence implies $d$-$\othree$ convergence.
\item Limits in $d$-$\othree$ are unique. Consequently, $d$-o limits and sequential limits are unique, and so are $d$-$\ru$ limits when $E$ is Archimedean.
\item Every convergent net (sequence) is Cauchy in the corresponding sense.
\item If $x_\lambda\to x$ and $y_\lambda\to y$ in $d$-$\othree$, with the same index set, then $d(x_\lambda,y_\lambda)\xrightarrow{\othree}d(x,y)$ in $E$.
\end{enumerate}
\end{proposition}

\begin{proof}
(i) The first two implications are immediate (take $\Gamma=\Lambda$). For $\ru$, index by $\eps\in(0,\infty)$ with the reverse order; since $E$ is Archimedean, $\inf_{\eps>0}\eps u=0$.

(ii) Let $x_\lambda\to x$ and $x_\lambda\to y$ with dominating nets $(a_\gamma)$ and $(b_\delta)$. Given $\gamma,\delta$, choose $\lambda$ beyond both thresholds; then $d(x,y)\le d(x,x_\lambda)+d(y,x_\lambda)\le a_\gamma+b_\delta$. If $v$ is a lower bound of $\{a_\gamma+b_\delta\}$, then $v-a_\gamma$ is a lower bound of $(b_\delta)$, so $v\le a_\gamma$ for all $\gamma$ and $v\le 0$. Hence $\inf_{(\gamma,\delta)}(a_\gamma+b_\delta)=0$ and $d(x,y)=0$.

(iii) Since $d(x_\lambda,x_\mu)\le d(x_\lambda,x)+d(x_\mu,x)$, use $2a_\gamma$, $2\eps u$, or $2a_n$ (for sequences: $d(x_n,x_{n+p})\le a_n+a_{n+p}\le 2a_n$).

(iv) By Lemma~\ref{a:lem:vmbasic}(iii), $\pm(d(x_\lambda,y_\lambda)-d(x,y))\le a_\gamma+b_\delta$ eventually, and $(a_\gamma+b_\delta)\downarrow 0$ over $\Gamma\times\Delta$ as in (ii).
\end{proof}

\begin{remark}\label{a:rem:lex}
Without the Archimedean property, $\ru$-limits need not be unique. Let $E=\R^2$ with the lexicographic cone and let $d(x,y)=(0,1)$ for $x\ne y$ on any set with at least two points. Then (VM1) and (VM2) hold, and $d(x,y)\le\eps(1,0)$ for every $\eps>0$, so every constant sequence $\ru$-converges to every point.
\end{remark}

\section{Transfer to vector lattice covers}\label{a:sec:transfer}

Throughout this section, $(Y,i)$ is a vector lattice cover of the pre-Riesz space $E$ and $(X,d,E)$ is a pre-Riesz-valued vector metric space.

\begin{proposition}[forward transfer]\label{a:prop:forward}
The map $i\id d$ is a $Y$-valued vector metric on $X$. Every net that converges (is Cauchy) in $d$-o, $d$-$\othree$ or $d$-$\ru$ converges (is Cauchy) in the corresponding sense for $i\id d$, with the same limit. The same holds for the sequential notions.
\end{proposition}

\begin{proof}
(VM1) follows from injectivity of $i$ and (VM2) from positivity and linearity. A dominating net $(a_\gamma)$ with $a_\gamma\downarrow 0$ in $E$ yields the dominating net $(i(a_\gamma))$, which decreases to $0$ in $Y$ by Lemma~\ref{a:lem:inf}. For $\ru$ replace $u$ by $i(u)$. For the order convergence of $E$ itself see also \cite[Proposition~3.7.8]{KvG19}.
\end{proof}

Relatively uniform convergence, in contrast, transfers in both directions without any hypothesis on the cover.

\begin{proposition}[$\ru$-invariance]\label{a:prop:ru}
For every cover $(Y,i)$, every net $(x_\lambda)$ in $X$ and every $x\in X$:
\begin{enumerate}[label=\textup{(\roman*)}]
\item $x_\lambda\to x$ in $d$-$\ru$ if and only if $x_\lambda\to x$ in $(i\id d)$-$\ru$;
\item $(x_\lambda)$ is $d$-$\ru$-Cauchy if and only if it is $(i\id d)$-$\ru$-Cauchy;
\item $X$ is $d$-$\ru$-complete if and only if it is $(i\id d)$-$\ru$-complete.
\end{enumerate}
\end{proposition}

\begin{proof}
The forward directions are Proposition~\ref{a:prop:forward}. Let $v\in Y_+$ and suppose that for every $\eps>0$ eventually $i(d(x_\lambda,x))\le\eps v$. By Lemma~\ref{a:lem:major} there is $u\in E$ with $v\le i(u)$; since $i(u)\ge 0$, $u\in E_+$. Then $i(d(x_\lambda,x))\le\eps\, i(u)$ and, by bipositivity, $d(x_\lambda,x)\le\eps u$ eventually. The Cauchy case is identical. Since $d$ and $i\id d$ live on the same set and have the same $\ru$-convergent and $\ru$-Cauchy nets, (iii) follows.
\end{proof}

The $\othree$ case requires replacing a dominating net in $Y$ by one in $E$. Under pervasiveness this is a result of Kalauch and van Gaans; under uniform density we give a direct construction.

\begin{lemma}[{\cite[Lemma~3.7.11]{KvG19}}]\label{a:lem:KvG}
Let $E$ be Archimedean and assume \textup{(P)}. If $(y_\gamma)$ is a net in $Y$ with $y_\gamma\downarrow 0$, then there is a net $(u_\mu)$ in $E$ with $u_\mu\downarrow 0$ such that for every $\mu$ there is $\gamma$ with $y_\gamma\le i(u_\mu)$.
\end{lemma}

For a vector lattice $E$ embedded as an order dense majorising sublattice of $Y$, the same domination appears in the proof of \cite[Corollary~3.2]{Bil23}. In that setting the forward passage of order convergence from $E$ to $Y$ is the regularity of the embedding, and regular sublattices are characterised in \cite[Section~3]{Bil23}; the reverse passage studied here is not treated there.

\begin{lemma}\label{a:lem:dense}
Let $E$ be Archimedean and assume \textup{(D)} with $u\in E_+$.
\begin{enumerate}[label=\textup{(\roman*)}]
\item If $(b_n)$ is a sequence in $Y$ with $b_n\downarrow 0$, then there is a sequence $(a_n)$ in $E$ with $a_n\downarrow 0$ and $b_n\le i(a_n)\le b_n+2^{1-n}i(u)$ for all $n$.
\item If $(y_\gamma)_{\gamma\in\Gamma}$ is a net in $Y$ with $y_\gamma\downarrow 0$, then there is a net $(a_\sigma)_{\sigma\in\Sigma}$ in $E$ with $a_\sigma\downarrow 0$ such that for every $\sigma$ there is $\gamma$ with $y_\gamma\le i(a_\sigma)$.
\end{enumerate}
\end{lemma}

\begin{proof}
By Lemma~\ref{a:lem:arch}, $Y$ is Archimedean. We first note: if $h\in E$, $y\in Y$ and $i(h)\le y+2^{1-n}i(u)$ for infinitely many $n$, then $i(h)\le y$. Indeed, for every $m$ there is such an $n$ with $2^{n}\ge m$, and then $m(i(h)-y)\le m2^{1-n}i(u)\le 2i(u)$; since $Y$ is Archimedean, $i(h)-y\le 0$.

(i) Apply (D) to $y=b_n+2^{-n}i(u)$ and $\eps=2^{-n}$ to obtain $a_n\in E$ with
\[
b_n+2^{-n}i(u)\le i(a_n)\le b_n+2^{1-n}i(u).
\]
Then $i(a_{n+1})\le b_{n+1}+2^{-n}i(u)\le b_n+2^{-n}i(u)\le i(a_n)$, so $(a_n)$ is decreasing, and $a_n\ge 0$. If $h\in E$ is a lower bound of $(a_n)$, then for $n\ge m$ we have $i(h)\le b_n+2^{1-n}i(u)\le b_m+2^{1-n}i(u)$. By the preliminary remark, $i(h)\le b_m$ for all $m$, so $i(h)\le 0$ and $h\le 0$. Thus $a_n\downarrow 0$.

(ii) Let $\Sigma=\Gamma\times\N$ with the relation $(\gamma,n)\preceq(\gamma',n')$ if either $(\gamma,n)=(\gamma',n')$, or $\gamma\le\gamma'$ and $n<n'$. This is a partial order, and $(\gamma_1,n_1),(\gamma_2,n_2)$ have the upper bound $(\gamma_3,\max\{n_1,n_2\}+1)$ for any $\gamma_3\ge\gamma_1,\gamma_2$; so $\Sigma$ is directed. As in (i) choose $a_{(\gamma,n)}\in E$ with
\[
y_\gamma+2^{-n}i(u)\le i(a_{(\gamma,n)})\le y_\gamma+2^{1-n}i(u).
\]
If $(\gamma,n)\prec(\gamma',n')$, then $n'\ge n+1$ and $i(a_{(\gamma',n')})\le y_{\gamma'}+2^{-n}i(u)\le y_\gamma+2^{-n}i(u)\le i(a_{(\gamma,n)})$. Hence the net is decreasing. A lower bound $h\in E$ satisfies $i(h)\le y_\gamma+2^{1-n}i(u)$ for all $n$, hence $i(h)\le y_\gamma$ for all $\gamma$, hence $h\le 0$. Finally, $y_\gamma\le i(a_{(\gamma,n)})$ by construction.
\end{proof}

We can now prove the transfer of order convergence of nets.

\begin{theorem}[$\othree$-invariance]\label{a:thm:o3}
Let $E$ be Archimedean and let $(Y,i)$ satisfy \textup{(P)} or \textup{(D)}. For every net $(x_\lambda)$ in $X$ and every $x\in X$:
\begin{enumerate}[label=\textup{(\roman*)}]
\item $x_\lambda\to x$ in $d$-$\othree$ if and only if $x_\lambda\to x$ in $(i\id d)$-$\othree$;
\item $(x_\lambda)$ is $d$-$\othree$-Cauchy if and only if it is $(i\id d)$-$\othree$-Cauchy;
\item $X$ is $d$-$\othree$-complete if and only if it is $(i\id d)$-$\othree$-complete.
\end{enumerate}
\end{theorem}

\begin{proof}
The forward directions are Proposition~\ref{a:prop:forward}. Under (P) the argument is that of \cite[Proposition~3.7.12]{KvG19}, applied to the distances. Let $(y_\gamma)\downarrow 0$ in $Y$ be such that for every $\gamma$ eventually $i(d(x_\lambda,x))\le y_\gamma$. Lemma~\ref{a:lem:KvG} under (P), or Lemma~\ref{a:lem:dense}(ii) under (D), gives a net $(a_\sigma)\downarrow 0$ in $E$ such that for each $\sigma$ there is $\gamma(\sigma)$ with $y_{\gamma(\sigma)}\le i(a_\sigma)$. Fix $\sigma$. Eventually $i(d(x_\lambda,x))\le y_{\gamma(\sigma)}\le i(a_\sigma)$, hence $d(x_\lambda,x)\le a_\sigma$ by bipositivity. This proves (i); (ii) is the same argument applied to $d(x_\lambda,x_\mu)$, and (iii) follows as in Proposition~\ref{a:prop:ru}.
\end{proof}

For sequences, uniform density suffices; pervasiveness does not, as Example~\ref{d:ex:counter} will show.

\begin{theorem}[sequential invariance]\label{a:thm:seq}
Let $E$ be Archimedean and let $(Y,i)$ satisfy \textup{(D)}. Then a sequence in $X$ converges (is Cauchy) in the sequential sense for $d$ if and only if it does so for $i\id d$. In particular, $X$ is complete in the sequential sense for $d$ if and only if it is for $i\id d$.
\end{theorem}

\begin{proof}
The forward direction is Proposition~\ref{a:prop:forward}. If $i(d(x_n,x))\le b_n$ with $b_n\downarrow 0$ in $Y$, then Lemma~\ref{a:lem:dense}(i) gives $a_n\downarrow 0$ in $E$ with $b_n\le i(a_n)$, so $d(x_n,x)\le a_n$. Cauchy sequences are treated in the same way.
\end{proof}

\begin{remark}\label{a:rem:quant}
Lemma~\ref{a:lem:dense}(i) is quantitative: the dominating sequence in $E$ exceeds the one in $Y$ by at most $2^{1-n}i(u)$. Hence also a priori estimates of the form $i(d(x_n,x))\le b_n$, as they appear in fixed point theorems, transfer to $d(x_n,x)\le a_n$ with an explicit error.
\end{remark}

As a consequence, results for vector-lattice-valued metrics can be imported.

\begin{corollary}\label{a:cor:import}
Let $E$ be Archimedean and let $(Y,i)$ satisfy \textup{(D)}. Then $(X,d)$ and $(X,i\id d)$ have the same convergent sequences, with the same limits, the same Cauchy sequences in the sequential sense of \cite{CA09}, and they are complete simultaneously. Consequently, every theorem on $Y$-valued vector metric spaces whose hypotheses and conclusions are formulated only in these terms (for instance the characterisations of continuity in \cite{Cev14}) applies to $(X,i\id d)$ and yields the same conclusion for $(X,d)$.
\end{corollary}

\section{Finite-dimensional value spaces}\label{a:sec:findim}

In finite dimensions the transfer holds for every cover. The proof uses the Kakutani--Krein--Yosida representation: an Archimedean vector lattice $Y$ with strong unit $e$ admits an injective lattice homomorphism $J\colon Y\to C(K)$, $K$ compact Hausdorff, with $J(e)=\mathbf 1$ and $J(Y)$ dense in $C(K)$ for the supremum norm; see \cite[Chapter~7]{LZ71}. %

\begin{lemma}\label{a:lem:dense-points}
Let $Y$, $e$, $J$, $K$ be as above and let $(y_\gamma)$ be a net in $Y$ with $y_\gamma\downarrow 0$. For every nonempty open set $W\subseteq K$ and every $\eps>0$ there are $\omega\in W$ and $\gamma$ with $J(y_\gamma)(\omega)<\eps$.
\end{lemma}

\begin{proof}
Suppose $J(y_\gamma)\ge\eps$ on $W$ for all $\gamma$. Choose $\omega_0\in W$ and, by Urysohn's lemma, $\varphi\in C(K)$ with $0\le\varphi\le\mathbf 1$, $\varphi(\omega_0)=1$ and $\varphi=0$ off $W$. By density pick $z\in Y$ with $\|J(z)-\varphi\|_\infty<\tfrac14$ and put $h=\bigl(z-\tfrac14 e\bigr)^+\wedge e\in Y$. Then $J(h)=(J(z)-\tfrac14)^+\wedge\mathbf 1$ vanishes off $W$, satisfies $0\le J(h)\le\mathbf 1$ and $J(h)(\omega_0)\ge\tfrac12$. Hence $\eps J(h)\le J(y_\gamma)$ on $K$ for all $\gamma$, so $\eps h$ is a lower bound of $(y_\gamma)$ with $\eps h>0$, a contradiction.
\end{proof}

The main result of this section shows that in finite dimensions no hypothesis on the cover is needed.

\begin{theorem}\label{a:thm:findim}
Let $E$ be finite-dimensional with a closed generating cone $E_+$, let $\|\cdot\|$ be a norm on $E$, and let $(Y,i)$ be any vector lattice cover. For a net $(e_\lambda)$ in $E_+$ the following are equivalent:
\begin{enumerate}[label=\textup{(\alph*)}]
\item $e_\lambda\xrightarrow{\ru}0$ in $E$;
\item $e_\lambda\xrightarrow{\othree}0$ in $E$;
\item $i(e_\lambda)\xrightarrow{\othree}0$ in $Y$;
\item $i(e_\lambda)\xrightarrow{\ru}0$ in $Y$;
\item $\|e_\lambda\|\to 0$.
\end{enumerate}
\end{theorem}

\begin{proof}
A closed cone makes $E$ Archimedean, and a generating cone in finite dimensions has an interior point $u$, which is an order unit; the order unit norm $\|\cdot\|_u$ is equivalent to $\|\cdot\|$. Hence (a)$\Leftrightarrow$(e), using that every $v\in E_+$ satisfies $v\le cu$ for some $c$. The implications (a)$\Rightarrow$(b) and (b)$\Rightarrow$(c) follow from Proposition~\ref{a:prop:relations}(i) and Lemma~\ref{a:lem:inf}, and (a)$\Leftrightarrow$(d) is the argument of Proposition~\ref{a:prop:ru} applied to elements. It remains to prove (c)$\Rightarrow$(e).

By Lemma~\ref{a:lem:major}, $e:=i(u)$ is a strong unit of $Y$, and $Y$ is Archimedean by Lemma~\ref{a:lem:arch}. Let $J\colon Y\to C(K)$ be as above. For $\omega\in K$ define $f_\omega\in E^*$ by $f_\omega(x)=J(i(x))(\omega)$. Each $f_\omega$ is positive with $f_\omega(u)=1$, and by bipositivity of $i$ and $J$,
\begin{equation}\label{a:eq:Kchar}
x\in E_+\iff f_\omega(x)\ge 0\ \text{ for all }\omega\in K.
\end{equation}
The map $\omega\mapsto f_\omega$ is continuous into $E^*$, since $E$ is finite-dimensional and each $\omega\mapsto f_\omega(x)$ is continuous.

Let $E_+^*$ be the dual cone and $B=\{f\in E_+^*:f(u)=1\}$; $B$ is compact because $|f(x)|\le\|x\|_u$ for $f\in B$. Let $F=\{f_\omega:\omega\in K\}\subseteq B$. We claim that $B=\overline{\mathrm{conv}}\,F$. Otherwise there are $f\in B$, $x\in E$ and $c\in\R$ with $f(x)<c\le g(x)$ for all $g\in\overline{\mathrm{conv}}\,F$. Then $g(x-cu)\ge 0$ for all $g\in F$, so $x-cu\in E_+$ by \eqref{a:eq:Kchar}, while $f(x-cu)<0$, contradicting $f\in E_+^*$.

Since $E_+$ is closed and $E_+\cap(-E_+)=\{0\}$, the dual cone $E_+^*$ has an interior point $f^*$, which we normalise by $f^*(u)=1$; then $f^*(x)\ge c_0\|x\|$ on $E_+$ for some $c_0>0$. Let $H=\{f\in E^*:f(u)=1\}$ and $V=\{f:f(u)=0\}$. As $f^*$ is interior to $E_+^*$, it lies in the relative interior of $B$ in $H$, and since the relative interior of a convex set and of its closure agree \cite[Theorem~6.3]{Roc70}, $f^*$ is in the relative interior of $\mathrm{conv}\,F$. Hence some ball $\{f\in H:\|f-f^*\|\le 2r\}$ is contained in $\mathrm{conv}\,F$. Choose a polytope in $H$ that contains the ball of radius $r$ about $f^*$ and is contained in the ball of radius $2r$; by Carath\'eodory's theorem each of its finitely many vertices is a convex combination of finitely many elements of $F$. So there are a finite set $F_0=\{f_{\omega_1},\dots,f_{\omega_m}\}\subseteq F$ and $r>0$ with $\{f\in H:\|f-f^*\|\le r\}\subseteq\mathrm{conv}\,F_0$.

\emph{Perturbation.} If $g_j'\in H$ satisfy $\|g_j'-f_{\omega_j}\|<r$ for all $j$, then $f^*\in\mathrm{conv}\{g_1',\dots,g_m'\}$. Otherwise, by separation in $H$ there is a linear functional $\ell$ on $E^*$, of norm one on $V$, with $\ell(g'_j)\le\ell(f^*)$ for all $j$. Pick $v\in V$ with $\|v\|\le 1$ and $\ell(v)=1$. Then $f^*+rv\in\mathrm{conv}\,F_0$, so $\ell(f^*)+r\le\max_j\ell(f_{\omega_j})\le\max_j\ell(g_j')+\max_j\|f_{\omega_j}-g_j'\|<\ell(f^*)+r$, a contradiction.

Now let $(y_\gamma)\downarrow 0$ in $Y$ dominate $i(e_\lambda)$ in the sense of $\othree$, and let $\eps>0$. Let $W_j=\{\omega:\|f_\omega-f_{\omega_j}\|<r\}$, an open neighbourhood of $\omega_j$. By Lemma~\ref{a:lem:dense-points} there are $\omega_j'\in W_j$ and $\gamma_j$ with $J(y_{\gamma_j})(\omega'_j)<\eps$. Choose $\gamma\ge\gamma_1,\dots,\gamma_m$ and write $f^*=\sum_j t_jf_{\omega'_j}$ with $t_j\ge 0$, $\sum_j t_j=1$. Eventually $i(e_\lambda)\le y_\gamma$, and then
\[
c_0\|e_\lambda\|\le f^*(e_\lambda)=\sum_j t_j\,J(i(e_\lambda))(\omega'_j)\le\sum_j t_j\,J(y_{\gamma_j})(\omega'_j)<\eps.
\]
\end{proof}

For vector metrics this reads as follows.

\begin{corollary}\label{a:cor:findim}
Let $E$ be finite-dimensional with a closed generating cone and let $(Y,i)$ be any vector lattice cover. For an $E$-valued vector metric $d$, a net converges (is Cauchy) in $d$-$\othree$, in $d$-$\ru$, in $(i\id d)$-$\othree$ or in $(i\id d)$-$\ru$ if and only if $\|d(x_\lambda,x)\|\to 0$ (respectively $\|d(x_\lambda,x_\mu)\|\to 0$). For sequences, the sequential notions of \cite{CA09} for $d$ and for $i\id d$ are equivalent to these as well.
\end{corollary}

\begin{proof}
Apply Theorem~\ref{a:thm:findim} to $e_\lambda=d(x_\lambda,x)$ and to $d(x_\lambda,x_\mu)$ (indexed by pairs). For sequences, if $\|d(x_n,x)\|_u\to 0$, then $a_n=\sup_{k\ge n}\|d(x_k,x)\|_u\,u$ is a decreasing sequence in $E$ with infimum $0$ dominating $d(x_n,x)$; the converse follows from Proposition~\ref{a:prop:relations}(i) and Theorem~\ref{a:thm:findim}.
\end{proof}

\section{Examples and sharpness}\label{a:sec:examples}

We begin with subspaces of $C(K)$, the basic source of uniformly dense covers.

\begin{proposition}\label{a:prop:CT}
Let $K$ be a compact Hausdorff space and $E\subseteq C(K)$ a subspace containing $\mathbf 1$ that is dense for the supremum norm, with the pointwise order. Then $E$ is an Archimedean pre-Riesz space, $(C(K),\mathrm{incl})$ is a vector lattice cover of $E$ satisfying \textup{(D)} with $u=\mathbf 1$, and for $x\in E$ the pair $\{x,-x\}$ has a supremum in $E$ only if $|x|\in E$. In particular, $E$ is a vector lattice only if it is a vector sublattice of $C(K)$.
\end{proposition}

\begin{proof}
The constant function $\mathbf 1$ is an order unit, so $E$ is directed, and it is Archimedean as a subspace of $C(K)$. Given $f\in C(K)$ and $\eps>0$, density gives $p\in E$ with $\|p-(f+2\eps)\|_\infty<\eps$, so $f+\eps\le p\le f+3\eps$. This yields (D), and it shows that $f$ is the infimum in $C(K)$ of $\{p\in E:p\ge f\}$, i.e.\ order density. If $s=\sup_E\{x,-x\}$ exists, then $s\ge|x|$, and $s\le p$ for all $p\in E$ with $p\ge|x|$; taking the pointwise infimum over such $p$ gives $s\le|x|$, so $s=|x|\in E$.
\end{proof}

Smooth functions satisfy both density conditions.

\begin{example}[$C^1{[0,1]}$: \textup{(P)} and \textup{(D)}]\label{a:ex:C1}
The space $E=C^1[0,1]$ satisfies the hypotheses of Proposition~\ref{a:prop:CT} with $K=[0,1]$. It is not a vector lattice: the function $x(t)=t-\tfrac12$ belongs to $E$, but $|x|$ is not differentiable at $t=\tfrac12$, so $|x|\notin E$, and by Proposition~\ref{a:prop:CT} the pair $\{x,-x\}$ has no supremum in $E$. The cover $C[0,1]$ also satisfies (P): if $f>0$ in $C[0,1]$, then $f\ge\eta>0$ on an interval, and a smooth bump of height $\eta$ supported there lies in $E$ between $0$ and $f$. Theorems~\ref{a:thm:o3} and \ref{a:thm:seq} apply.
\end{example}

Polynomials satisfy \textup{(D)} but not \textup{(P)}.

\begin{example}[polynomials: \textup{(D)} without \textup{(P)}]\label{a:ex:poly}
Let $E=\mathcal P[0,1]$ be the polynomials on $[0,1]$. By the Weierstrass theorem and Proposition~\ref{a:prop:CT}, $(C[0,1],\mathrm{incl})$ is a cover satisfying (D), and $E$ is not a vector lattice. It does not satisfy (P): a continuous bump $f>0$ vanishing on an interval admits no polynomial $p$ with $0<p\le f$, since such $p$ would vanish on an interval. Theorems~\ref{a:thm:o3} and \ref{a:thm:seq} apply, but Lemma~\ref{a:lem:KvG} does not.
\end{example}

In finite dimensions both conditions may fail, although the transfer holds by Corollary~\ref{a:cor:findim}.

\begin{example}[a pyramid cone: neither \textup{(P)} nor \textup{(D)}]\label{a:ex:pyramid}
Let $E=\R^3$ with $E_+=\{(a,b,c):c\ge|a|,\ c\ge|b|\}$, a closed generating cone with the four extreme rays spanned by $(\pm1,\pm1,1)$. Since a finite-dimensional Archimedean vector lattice is lattice isomorphic to $\R^n$ with the standard cone, which has exactly $n$ extreme rays, $E$ is not a vector lattice. Define $i\colon E\to\R^4$ by $i(a,b,c)=(c-a,\,c+a,\,c-b,\,c+b)$. Then $i$ is injective and bipositive by the definition of $E_+$, and $i(E)=\{y:y_1+y_2=y_3+y_4\}$. For $y\in\R^4$ with $\delta=(y_1+y_2)-(y_3+y_4)\ge 0$, the vectors $y+\delta e_3$ and $y+\delta e_4$ lie in $i(E)$, lie above $y$, and have infimum $y$; the case $\delta<0$ is symmetric. So $(\R^4,i)$ is a vector lattice cover. It does not satisfy (P): if $0\le i(x)\le e_1$, then $c+a=c-b=c+b=0$, which forces $x=0$. It does not satisfy (D) because $i(E)$ is a proper closed subspace. Nevertheless, by Theorem~\ref{a:thm:findim} the conclusions of Theorems~\ref{a:thm:o3} and \ref{a:thm:seq} hold. Thus neither (P) nor (D) is necessary.
\end{example}

We now show that the transfer may fail when neither \textup{(P)} nor \textup{(D)} holds. Let $Y$ be the space of bounded real sequences $y=(y_j)_{j\in\Z}$ for which $\lim_{j\to+\infty}y_j$ exists, with the pointwise order; $Y$ is an Archimedean vector lattice. Let
\[
\varphi(y)=\lim_{j\to+\infty}y_j-\sum_{k=1}^\infty 2^{-k}y_{-k},\qquad Z=\ker\varphi,
\]
with the order inherited from $Y$, and $i\colon Z\to Y$ the inclusion. Write $e^{(j)}$ for the $j$-th unit sequence. This is \cite[Example~3.7.10]{KvG19}, where it is shown that $(e^{(m)})$ order converges to $0$ in $Y$ but not in $Z$. Theorem~\ref{a:thm:counter} below strengthens this: the sequence is not even $\othree$-Cauchy for the vector metric $d$.

\begin{lemma}\label{a:lem:Z}
The pair $(Y,i)$ is a vector lattice cover of the Archimedean pre-Riesz space $Z$, and it satisfies neither \textup{(P)} nor \textup{(D)}.
\end{lemma}

\begin{proof}
Bipositivity is clear and $Z$ is Archimedean. Let $y\in Y$ and $\alpha=\varphi(y)$. If $\alpha\ge 0$, then $x_1=y+2\alpha e^{(-1)}$ and $x_2=y+4\alpha e^{(-2)}$ lie in $Z$, lie above $y$, and $x_1\wedge x_2=y$. If $\alpha<0$, then $x_m=y+|\alpha|\sum_{j\ge m}e^{(j)}$ lies in $Z$, lies above $y$, and $\inf_m x_m=y$. Hence $Z$ is order dense in $Y$, and $Z$ is pre-Riesz by van Haandel's theorem. If $z\in Z$ and $0\le z\le e^{(-1)}$, then $z=te^{(-1)}$ and $0=\varphi(z)=-t/2$, so (P) fails. Since $|\varphi(y)|\le 2\|y\|_\infty$, the kernel $Z$ is a proper closed subspace for the norm of the strong unit $\mathbf 1$, so (D) fails.
\end{proof}

For $m\in\N$ we have $e^{(m)}\in Z$. Let $X=\{0\}\cup\{e^{(m)}:m\in\N\}\subseteq Z_+$ and
\[
d(x,y)=\begin{cases}x+y,&x\ne y,\\ 0,&x=y.\end{cases}
\]
(VM1) holds because distinct points of $X$ have nonzero sum. For (VM2) with $x\neq y$: if $z\notin\{x,y\}$, then $d(x,z)+d(y,z)=x+y+2z\ge x+y$, and if $z\in\{x,y\}$ the right-hand side equals $x+y$.

\begin{theorem}\label{a:thm:counter}
In the $Z$-valued vector metric space $(X,d,Z)$:
\begin{enumerate}[label=\textup{(\roman*)}]
\item $e^{(m)}\to 0$ for $i\id d$ in the sequential sense, hence in $(i\id d)$-o and $(i\id d)$-$\othree$;
\item $(e^{(m)})$ is not $d$-$\othree$-Cauchy, and in particular it does not $d$-$\othree$-converge.
\end{enumerate}
\end{theorem}

\begin{proof}
(i) Let $s_m=\sum_{j\ge m}e^{(j)}\in Y$. Then $s_m\downarrow 0$ in $Y$ and $d(e^{(m)},0)=e^{(m)}\le s_m$.

(ii) Suppose $(w_\gamma)\downarrow 0$ in $Z$ and for each $\gamma$ there is $m_\gamma$ with $e^{(m)}+e^{(n)}\le w_\gamma$ for all distinct $m,n\ge m_\gamma$. Then $w_\gamma\ge 0$ and $(w_\gamma)_j\ge 1$ for all $j\ge m_\gamma$, so $\lim_{j\to+\infty}(w_\gamma)_j\ge 1$, and $\varphi(w_\gamma)=0$ gives $\sum_k2^{-k}(w_\gamma)_{-k}\ge 1$. Fix $\gamma_0$ and let $c_k=\inf_{\gamma\ge\gamma_0}(w_\gamma)_{-k}\ge 0$. For every $N$ and $\gamma\ge\gamma_0$,
\[
\sum_{k=1}^N2^{-k}(w_\gamma)_{-k}\ge 1-2^{-N}\|w_{\gamma_0}\|_\infty.
\]
The left-hand side is a finite sum of decreasing real nets, so the same bound holds for $\sum_{k\le N}2^{-k}c_k$. Letting $N\to\infty$ gives $\sum_k2^{-k}c_k\ge 1$, so $\delta:=c_k>0$ for some $k$. Put $w=\delta e^{(-k)}-2\delta e^{(-k-1)}$. Then $\varphi(w)=-(\delta2^{-k}-2\delta2^{-k-1})=0$, so $w\in Z$, and $w\le w_\gamma$ for all $\gamma\ge\gamma_0$, hence for all $\gamma$. Since $w\not\le 0$, this contradicts $\inf_\gamma w_\gamma=0$ in $Z$. The last assertion follows from Proposition~\ref{a:prop:relations}(iii).
\end{proof}

\begin{remark}\label{a:rem:counter}
In Theorem~\ref{a:thm:counter} even a dominating \emph{sequence} in $Y$ cannot be replaced by any dominating \emph{net} in $Z$. The failure is therefore caused by the cover, not by the distinction between sequences and nets.
\end{remark}

\begin{table}[ht]
\centering
\caption{Covers and the $\othree$-transfer}\label{a:tab:summary}
\small
\begin{tabular}{@{}lcccc@{}}
\toprule
Value space $E$ and cover & (P) & (D) & $\othree$-transfer & sequential transfer\\
\midrule
$C^1[0,1]\subseteq C[0,1]$ & yes & yes & yes & yes\\
$\mathcal P[0,1]\subseteq C[0,1]$ & no & yes & yes & yes\\
pyramid cone, $\R^3\to\R^4$ & no & no & yes & yes\\
$Z\subseteq Y$ (Lemma~\ref{a:lem:Z}) & no & no & no & no\\
$\R\mathbf 1\oplus c_{00}(\Gamma)\subseteq\ell^\infty(\Gamma)$, $\Gamma$ uncountable (Ex.~\ref{d:ex:counter}) & yes & no & yes & no\\
$\R\mathbf 1\oplus c_{00}(\N)\subseteq\ell^\infty$ (Ex.~\ref{d:ex:counter}) & yes & no & yes & yes\\
\bottomrule
\end{tabular}
\end{table}

\section{Nets: a characterisation of the transfer}\label{d:sec:nets}

Throughout Sections~\ref{d:sec:nets} and~\ref{d:sec:seq}, $E$ is an Archimedean pre-Riesz space with cone $E_+$ and vector lattice cover $(Y,i)$; $Y$ is then Archimedean (Lemma~\ref{a:lem:arch}). A net $(e_\lambda)$ in $E_+$ is \emph{$\othree$-null}, $e_\lambda\ocv 0$, if $e_\lambda\xrightarrow{\othree}0$ in the sense of Definition~\ref{a:def:conv}, and a sequence $(e_n)$ in $E_+$ is \emph{order null} if $e_n\le a_n$ for some sequence $a_n\downarrow 0$; thus $x_\lambda\to x$ in $d$-$\othree$ if and only if $d(x_\lambda,x)\ocv 0$, and $x_n\to x$ in the sequential sense if and only if $(d(x_n,x))$ is order null (Definition~\ref{a:def:vmconv}).

\begin{lemma}[star spaces]\label{d:lem:star}
For a family $(h_s)_{s\in S}$ in $E_+$, the set $X_\star=\{\xi\}\cup\{p_s:h_s\ne 0\}$ with $d(p_s,\xi)=h_s$ and $d(p_s,p_t)=h_s+h_t$ ($s\ne t$) is an $E$-valued vector metric space. With $\pi(s)=p_s$ if $h_s\ne 0$ and $\pi(s)=\xi$ otherwise, a net or sequence $\pi(s_\lambda)$ converges to $\xi$ for $d$, respectively for $i\id d$, if and only if $(h_{s_\lambda})$ is null in $E$, respectively $(i(h_{s_\lambda}))$ is null in $Y$, in the corresponding sense.
\end{lemma}

\begin{proof}
For (VM2) with distinct $P,Q$: if $R=\xi$, the inequality is an equality, and if $R=p_r$ the right-hand side exceeds the left by $2h_r\ge 0$. The rest is immediate from $d(\pi(s),\xi)=h_s$.
\end{proof}

The transfer of order convergence is governed by two domination conditions on the cover.

\begin{definition}\label{d:def:CD}
The cover $(Y,i)$ has
\begin{enumerate}
\item[\CD] if for every net $y_\gamma\downarrow 0$ in $Y$ there is a net $a_\sigma\downarrow 0$ in $E$ such that for every $\sigma$ there is $\gamma$ with $y_\gamma\le i(a_\sigma)$;
\item[\CDf] if for every net $y_\gamma\downarrow 0$ in $Y$ there is a net $a_\sigma\downarrow 0$ in $E$ such that for every $\sigma$ there is $\gamma$ with $e\le a_\sigma$ for all $e\in E_+$ with $i(e)\le y_\gamma$.
\end{enumerate}
\end{definition}

The second condition characterises the transfer of order convergence of nets.

\begin{theorem}\label{d:thm:charact}
The following are equivalent:
\begin{enumerate}[label=\textup{(\alph*)}]
\item for every $E$-valued vector metric space, $d$-$\othree$ and $(i\id d)$-$\othree$ convergence of nets coincide, and so do the corresponding Cauchy nets;
\item for every net $(e_\lambda)$ in $E_+$, $i(e_\lambda)\ocv 0$ in $Y$ implies $e_\lambda\ocv 0$ in $E$;
\item the cover has \CDf.
\end{enumerate}
\end{theorem}

\begin{proof}
(a)$\Rightarrow$(b). Apply (a) to the star space of Lemma~\ref{d:lem:star} with $S=\Lambda$ and $h_\lambda=e_\lambda$.

(b)$\Rightarrow$(c). Let $y_\gamma\downarrow 0$ in $Y$. On $\Lambda=\{(\gamma,e):e\in E_+,\ i(e)\le y_\gamma\}$ put $(\gamma,e)\preceq(\gamma',e')$ if $\gamma\le\gamma'$; this is a directed preorder. The net $e_{(\gamma,e)}=e$ satisfies $i(e_{(\gamma',e')})\le y_{\gamma'}\le y_\gamma$ whenever $(\gamma',e')\succeq(\gamma,e)$, so $i(e_\lambda)\ocv 0$. By (b) there is $a_\sigma\downarrow 0$ in $E$ such that for each $\sigma$ there is $(\gamma_\sigma,e_\sigma)$ with $e\le a_\sigma$ for all $(\gamma,e)\succeq(\gamma_\sigma,e_\sigma)$. In particular, $e\le a_\sigma$ for every $e\in E_+$ with $i(e)\le y_{\gamma_\sigma}$.

(c)$\Rightarrow$(a). The forward implication always holds (Proposition~\ref{a:prop:forward}). If $i(d(x_\lambda,x))\le y_\gamma$ eventually, for every $\gamma$, with $y_\gamma\downarrow 0$, take $a_\sigma$ from \CDf. For each $\sigma$ there is $\gamma$ such that eventually $e=d(x_\lambda,x)\in E_+$ satisfies $i(e)\le y_\gamma$, hence $d(x_\lambda,x)\le a_\sigma$. Cauchy nets are treated in the same way with $d(x_\lambda,x_\mu)$.
\end{proof}

All sufficient conditions obtained so far are covered.

\begin{corollary}\label{d:cor:sufficient}
Condition \CDf\ holds in each of the following cases: $i(E)$ is pervasive; $i(E)$ is uniformly dense; $\dim E<\infty$ and $E_+$ is closed and generating. It fails for the cover of Theorem~\ref{a:thm:counter}.
\end{corollary}

\begin{proof}
Under (P) or (D), Lemmas~\ref{a:lem:KvG} and~\ref{a:lem:dense} give \CD, which implies \CDf\ by bipositivity: $i(e)\le y_\gamma\le i(a_\sigma)$ gives $e\le a_\sigma$. In finite dimensions, Theorem~\ref{a:thm:findim} gives (b) of Theorem~\ref{d:thm:charact}. The last statement follows from Theorem~\ref{d:thm:charact}, since the space of Theorem~\ref{a:thm:counter} violates (a).
\end{proof}

The two conditions are related as follows.

\begin{proposition}\label{d:prop:CD}
Condition \CD\ implies \CDf. The converse holds if the cover has \emph{positive inner density}: $y=\sup\{i(e):e\in E_+,\ i(e)\le y\}$ for every $y\in Y_+$.
\end{proposition}

\begin{proof}
The first statement was shown in the proof of Corollary~\ref{d:cor:sufficient}. Under positive inner density, the $a_\sigma$ of \CDf\ satisfy $i(a_\sigma)\ge\sup\{i(e):e\in E_+,\ i(e)\le y_\gamma\}=y_\gamma$.
\end{proof}

\begin{remark}\label{d:rem:PID}
Positive inner density can fail even when \CD\ holds. For the cover $\R^3\to\R^4$ of the pyramid cone in Example~\ref{a:ex:pyramid}, no $e\in E_+\setminus\{0\}$ satisfies $i(e)\le e_1$, so the supremum in Proposition~\ref{d:prop:CD} is $0$ for $y=e_1$. On the other hand, \CD\ holds there, since $i(E)$ contains the strong unit $(1,1,1,1)$ of $\R^4$ and nets decreasing to $0$ in $\R^4$ converge in norm. Whether \CDf\ and \CD\ are always equivalent remains open (Problem~\ref{d:prob:CD}).
\end{remark}

\section{Sequences: countable domination and a pervasive counterexample}\label{d:sec:seq}

For sequences we use the sequential analogue of \CD.

\begin{definition}\label{d:def:CDs}
The cover has \CDs\ if for every sequence $b_n\downarrow 0$ in $Y$ there is a sequence $a_n\downarrow 0$ in $E$ with $b_n\le i(a_n)$ for all $n$.
\end{definition}

For sequences, the corresponding condition is sufficient, and the transfer can again be characterised by star spaces.

\begin{theorem}\label{d:thm:seq}
\begin{enumerate}[label=\textup{(\roman*)}]
\item Under \CDs, sequential convergence and Cauchy sequences for $d$ and for $i\id d$ coincide in every $E$-valued vector metric space, and so does completeness.
\item Uniform density implies \CDs.
\item If $i(E)$ contains a strong unit $u$ of $Y$ and every sequence $b_n\downarrow 0$ in $Y$ is $u$-uniformly null, that is, $b_n\le\eps_nu$ for some reals $\eps_n\downarrow 0$, then \CDs\ holds. This applies when $\dim Y<\infty$ and $E_+$ has an interior point $e_0$: then $i(e_0)$ is a strong unit of $Y$, because $i(E)$ majorises $Y$, and sequences decreasing to $0$ in $Y$ converge in norm.
\item Sequential convergence transfers for all $E$-valued vector metrics if and only if every sequence $(e_n)$ in $E_+$ with $(i(e_n))$ order null in $Y$ is order null in $E$.
\end{enumerate}
\end{theorem}

\begin{proof}
(i) If $i(d(x_n,x))\le b_n\downarrow 0$, take $a_n$ from \CDs; then $i(d(x_n,x))\le i(a_n)$, so $d(x_n,x)\le a_n$. (ii) is Lemma~\ref{a:lem:dense}(i). (iii) If $b_n\le\eps_nu$ with $\eps_n\downarrow 0$, put $a_n=\eps_ni^{-1}(u)$. (iv) follows from Lemma~\ref{d:lem:star}.
\end{proof}

The following example separates nets from sequences.

\begin{example}[a pervasive cover without sequential transfer]\label{d:ex:counter}
Let $\Gamma$ be an infinite set, $Y=\ell^\infty(\Gamma)$, and $E=\R\mathbf 1\oplus c_{00}(\Gamma)$, the functions that are constant off a finite set, with $i$ the inclusion.
\begin{enumerate}[label=(\alph*)]
\item The space $E$ is an Archimedean vector lattice with strong unit $\mathbf 1$, and $(Y,i)$ is a vector lattice cover. For $y\in Y$ and $t\in\Gamma$, the function $x_t=\|y\|_\infty\mathbf 1+(y(t)-\|y\|_\infty)e_t\in E$ satisfies $x_t\ge y$ and $x_t(t)=y(t)$, so the infimum of the elements of $E$ above $y$ is the pointwise infimum, $y$.
\item (P) holds: if $y>0$, then $y(t)>0$ for some $t$, and $0<y(t)e_t\le y$. (D) fails: for an infinite, coinfinite $A\subseteq\Gamma$ and $y=\mathbf 1_A$, any $x\in E$ with $y\le x\le y+\eps u$ has constant part at least $1$, so $u\ge 1/\eps$ off a finite subset of $\Gamma\setminus A$; this cannot hold for all $\eps>0$.
\item Condition \CD\ holds directly. Given $y_\gamma\downarrow 0$, index by pairs $(\gamma,F)$ with $F\subseteq\Gamma$ finite, ordered componentwise, and put $a_{(\gamma,F)}=\|y_\gamma\|_\infty\mathbf 1_{\Gamma\setminus F}+\sum_{t\in F}y_\gamma(t)e_t\in E$. This net decreases, dominates $y_\gamma$, and has pointwise infimum $0$. Hence $\othree$-convergence of nets transfers by Theorem~\ref{d:thm:charact}.
\item Let $\Gamma$ be uncountable and $t_1,t_2,\dots$ distinct points of $\Gamma$. Put $h_n=e_{t_n}$ and $b_n=\mathbf 1_{\{t_m:m\ge n\}}\in Y$. Then $h_n\le b_n\downarrow 0$, so in the star space of Lemma~\ref{d:lem:star} we have $p_n\to\xi$ for $i\id d$ in the sequential sense. Suppose $h_n\le a_n$ with $a_n\downarrow 0$ in $E$. Since $(a_n)$ decreases, $a_n\ge h_m$ for $m\ge n$, so $a_n(t_m)\ge 1$ for infinitely many $m$. Writing $a_n=\kappa_n\mathbf 1+g_n$ with $g_n$ finitely supported gives $\kappa_n\ge 1$. Let $S=\bigcup_n\operatorname{supp}g_n$, a countable set, and pick $t\in\Gamma\setminus S$. Then $a_n(t)=\kappa_n\ge 1$ for all $n$, so $e_t\in E$ is a positive lower bound of $(a_n)$, a contradiction. Hence $p_n\not\to\xi$ in the sequential sense for $d$.
\item Let $\Gamma=\N$. Then \CDs\ holds, and sequential convergence transfers although (D) fails. Given $b_n\downarrow 0$ in $\ell^\infty$, let $\kappa_n=\sup_{j\ge n}b_n(j)$, $a_n(j)=b_n(j)$ for $j<n$ and $a_n(j)=\kappa_n$ for $j\ge n$. Then $a_n\in E$, $a_n\ge b_n$, $(a_n)$ decreases, and $a_n(j)=b_n(j)\to 0$ for each $j$.
\end{enumerate}
\end{example}

We summarise the consequences of the example.

\begin{corollary}\label{d:cor:answer}
Pervasiveness implies the transfer of $\othree$-convergence of nets, but not of sequential convergence. Moreover, the dominating net in \textup{\CD} cannot in general be chosen with the index set of the given net: in Example~\ref{d:ex:counter}(d) with $\Gamma$ uncountable, the sequence $(h_n)$ is dominated by $b_n\downarrow 0$ in $Y$, but by no sequence decreasing to $0$ in $E$.
\end{corollary}

In contrast, uniform density does suffice for sequences (Theorem~\ref{a:thm:seq}).

\begin{remark}\label{d:rem:dedekind}
In Example~\ref{d:ex:counter}, $E$ is a vector lattice and $Y$ is its Dedekind completion, since $\mathbf 1\in E$ and $E$ is order dense in the Dedekind complete space $Y$. Thus sequential order convergence need not pass from the Dedekind completion of a vector lattice back to the vector lattice, although $\othree$-convergence of nets does. The same pair shows that the diagonal property of order convergence need not pass from a vector lattice to its Dedekind completion \cite[Example~7.6]{DC}.

The dichotomy between countable and uncountable $\Gamma$ is governed by order separability. By \cite[Theorem~29.3]{LZ71}, an Archimedean vector lattice is order separable if and only if, for every $u\ge 0$ and every disjoint family of nonzero positive elements, at most countably many members of the family are not disjoint from $u$; using the family $(e_t)_{t\in\Gamma}$ and $u=\mathbf 1$, one sees that $E$ is order separable if and only if $\Gamma$ is countable. For an order separable Archimedean vector lattice $L$ with Dedekind completion $L^\delta$, the proof of \cite[Theorem~32.9]{LZ71} shows that for every sequence $0\le w_n\uparrow w_0$ in $L^\delta$ there are $u_j\in L$ with $0\le u_j\uparrow w_0$ in $L^\delta$ and indices $k(j)$ with $u_j\le w_{k(j)}$; then $v_n=\sup\{u_j:j\le n,\ k(j)\le n\}$ (and $v_n=0$ if this set is empty) satisfies $0\le v_n\uparrow w_0$ and $v_n\le w_n$ for all $n$; see also \cite[Exercise~75.7]{LZ71}. Given $b_n\downarrow 0$ in $L^\delta$, choose $w_0\in L$ with $w_0\ge b_1$ and apply this to $w_n=w_0-b_n$; then $a_n=w_0-v_n\in L$ satisfies $a_n\downarrow 0$ and $b_n\le a_n$. Hence \CDs\ holds for every order separable vector lattice in its Dedekind completion, which explains Example~\ref{d:ex:counter}(e).
\end{remark}

\section{Integral vector metrics}\label{d:sec:integral}

Let $E=C^1[0,1]$ with the pointwise order and $Y=C[0,1]$; the pair $(Y,\mathrm{incl})$ is a cover satisfying (P) and (D) (Example~\ref{a:ex:C1}). Vector metrics defined by integral kernels are the simplest $E$-valued metrics that do not take values in a vector lattice in any natural way.

\begin{example}\label{a:ex:integral}
Let $k\colon[0,1]^2\to[0,\infty)$ be continuous with $\partial k/\partial t$ continuous and $k(s,s)>0$ for all $s$. For $f,g\in X=C[0,1]$ put
\[
d(f,g)(t)=\int_0^1k(t,s)\,|f(s)-g(s)|\,ds.
\]
Then $d(f,g)\in C^1[0,1]$, (VM2) holds pointwise, and (VM1) holds since $d(f,g)=0$ forces $f=g$ near each $s$. So $d$ is a $C^1[0,1]$-valued vector metric.
\end{example}

More generally, for a measurable kernel $k\colon[0,1]^2\to[0,\infty)$ such that $t\mapsto\int_0^1k(t,s)\varphi(s)\,ds$ belongs to $C^1[0,1]$ for every $\varphi\in C[0,1]$, we write $d_k(f,g)(t)=\int_0^1k(t,s)|f(s)-g(s)|\,ds$; this is a vector metric on $C[0,1]$ whenever $d_k(f,g)=0$ forces $f=g$. Since $|f-g|\in C[0,1]$, the hypothesis on $k$ guarantees $d_k(f,g)\in C^1[0,1]$. The Volterra kernel $k(t,s)=(t-s)_+$ of Theorem~\ref{d:thm:volterra} is an instance: it is only piecewise smooth across the diagonal, but $\frac{d}{dt}\int_0^1(t-s)_+\varphi(s)\,ds=\int_0^t\varphi(s)\,ds$ is continuous.

\begin{lemma}\label{d:lem:dense-zero}
If $(w_\gamma)$ is a net in $C^1[0,1]$ with $w_\gamma\downarrow 0$, then $g=\inf_\gamma w_\gamma$ (pointwise) vanishes on a dense subset of $[0,1]$.
\end{lemma}

\begin{proof}
By Lemma~\ref{a:lem:inf}, $w_\gamma\downarrow 0$ in $C[0,1]$. If $g\ge\eps$ on an open interval $I$, then a continuous bump of height $\eps$ supported in $I$ is a positive lower bound. So each $\{g<\eps\}$ is dense; it is open since $g$ is upper semicontinuous, and Baire's theorem gives density of $\{g=0\}$.
\end{proof}

If the kernel is positive on the diagonal, nothing is gained from the vector-valued structure.

\begin{theorem}[positive diagonal: no degeneracy]\label{d:thm:band}
Let $k$ be continuous, nonnegative, with continuous $\partial k/\partial t$ and $k(s,s)>0$ for all $s$; no assumption is made on $k$ away from the diagonal, so $k$ may vanish outside a narrow band $\{|t-s|\le\delta\}$. For a net $(f_\lambda)$ and $f$ in $X$ the following are equivalent: $f_\lambda\to f$ in $d_k$-$\ru$; $f_\lambda\to f$ in $d_k$-$\othree$; $\|f_\lambda-f\|_{L^1}\to 0$. For sequences, these are also equivalent to convergence in the sequential sense.
\end{theorem}

\begin{proof}
By compactness, there are $c,\delta>0$ with $k(t,s)\ge c$ whenever $|t-s|\le\delta$; let $M=\max k$. If $c_\lambda=\|f_\lambda-f\|_{L^1}\to 0$, then $d_k(f_\lambda,f)\le Mc_\lambda\mathbf 1$, which gives $\ru$-convergence, hence $\othree$-convergence (Proposition~\ref{a:prop:relations}), and for sequences convergence in the sequential sense with $a_n=M\sup_{j\ge n}c_j\mathbf 1$.

Conversely, let $(w_\gamma)\downarrow 0$ dominate $d_k(f_\lambda,f)$ in the sense of $\othree$, and let $\eps>0$. By Lemma~\ref{d:lem:dense-zero}, applied to $g=\inf_\gamma w_\gamma$, there are points $t_1,\dots,t_N$ with $g(t_j)=0$ such that the intervals $[t_j-\delta,t_j+\delta]$ cover $[0,1]$; here $N$ depends only on $\delta$. Choose $\gamma_j$ with $w_{\gamma_j}(t_j)<\eps$ and $\gamma\ge\gamma_1,\dots,\gamma_N$. Eventually $d_k(f_\lambda,f)\le w_\gamma$, and then
\[
c\int_{|s-t_j|\le\delta}|f_\lambda-f|\,ds\le d_k(f_\lambda,f)(t_j)<\eps\qquad(j=1,\dots,N),
\]
so $\|f_\lambda-f\|_{L^1}\le N\eps/c$. Convergence in the sequential sense implies $\othree$-convergence.
\end{proof}

Thus kernels supported in a band around the diagonal but positive on it give nothing beyond $L^1$; the case $0<m\le k\le M$, in which $d_k$ is squeezed between multiples of $\mathbf 1$, is a special case. Degeneracy therefore requires a kernel that vanishes somewhere on the diagonal; the model case is the Volterra kernel.

\begin{theorem}[Volterra kernel]\label{d:thm:volterra}
Let $k(t,s)=(t-s)_+$, so that $d(f,g)(t)=\int_0^t(t-s)|f(s)-g(s)|\,ds$ is a $C^1[0,1]$-valued vector metric on $X=C[0,1]$. For a sequence $(f_n)$ and $f$ in $X$, write $\nu(h)=\int_0^1(1-s)|h(s)|\,ds$.
\begin{enumerate}[label=\textup{(\alph*)}]
\item $f_n\to f$ in $d$-$\ru$ if and only if $\nu(f_n-f)\to 0$.
\item $f_n\to f$ in the sequential sense if and only if $\int_0^t|f_n-f|\to 0$ for every $t<1$ and $\sup_n\nu(f_n-f)<\infty$.
\item $L^1$-convergence implies (a), which implies (b), and neither implication can be reversed.
\end{enumerate}
\end{theorem}

\begin{proof}
The function $d(f,g)$ is nondecreasing and convex in $t$, with $d(f,g)'(t)=\int_0^t|f-g|$; it vanishes identically only if $f=g$, and $\max_td(f,g)=d(f,g)(1)=\nu(f-g)$. Write $D_n=d(f_n,f)$.

(a) If $D_n\le\eps_nu$ with $\eps_n\downarrow 0$, then $\nu(f_n-f)=D_n(1)\le\eps_nu(1)\to 0$. Conversely, $D_n\le\nu(f_n-f)\mathbf 1$.

(b) Necessity. Let $D_n\le a_n\downarrow 0$ in $C^1[0,1]$. Then $\nu(f_n-f)=D_n(1)\le a_1(1)$. By Lemma~\ref{d:lem:dense-zero}, $\inf_na_n(t')=0$ on a dense set of $t'$. For $t<t'<1$ with $t'$ in that set,
\[
\int_0^t|f_n-f|\le\frac{D_n(t')}{t'-t}\le\frac{a_n(t')}{t'-t}\to 0 .
\]
Sufficiency. Let $M=\max\{1,\sup_n\nu(f_n-f)\}$ and $G_n=\sup_{m\ge n}D_m$ (pointwise). Then $G_n\le M$ is nondecreasing in $t$, decreasing in $n$, and $G_n(t)\to 0$ for every $t<1$, since $D_m(t)\le t\int_0^t|f_m-f|$. Choose $N_1<N_2<\cdots$ with $G_{N_j}(1-2^{-j})\le 2^{-j}$. Fix smooth $\chi_j\colon[0,1]\to[0,1]$ with $\chi_j=0$ on $[0,1-2^{1-j}]$ and $\chi_j=1$ on $[1-2^{-j},1]$; then $\chi_{j+1}\le\chi_j$. Put $a_n=2M\mathbf 1$ for $n<N_1$ and $a_n=2^{-j}\mathbf 1+M\chi_j$ for $N_j\le n<N_{j+1}$. Each $a_n\in C^1[0,1]$, and $(a_n)$ decreases. For $N_j\le n<N_{j+1}$: on $[0,1-2^{-j}]$, $D_n\le G_n\le G_{N_j}(1-2^{-j})\le 2^{-j}$; on $[1-2^{-j},1]$, $D_n\le M=M\chi_j$. So $D_n\le a_n$. For $t<1$, $\chi_j(t)=0$ for large $j$, hence $a_n(t)\to 0$. So the pointwise infimum vanishes on $[0,1)$, and $a_n\downarrow 0$ in $C[0,1]$, hence in $C^1[0,1]$ (Lemma~\ref{a:lem:inf}).

(c) The inequalities $\nu(h)\le\|h\|_{L^1}$ and $\int_0^t|h|\le(1-t)^{-1}\nu(h)$ give the implications. Let $\varphi$ be the tent on $[0,1]$ with $\varphi(0)=1$, $\varphi(1)=0$ extended by $0$, and $\varphi_n(s)=\varphi(n(1-s))$, supported in $[1-1/n,1]$. Then $\|n\varphi_n\|_{L^1}=\tfrac12$ while $\nu(n\varphi_n)=O(1/n)$: (a) holds for $f_n=n\varphi_n\to 0$ without $L^1$-convergence. And $\nu(n^2\varphi_n)=\int_0^1r\varphi(r)\,dr>0$ is constant while $\int_0^tn^2\varphi_n=0$ for $n>1/(1-t)$: (b) holds for $f_n=n^2\varphi_n\to 0$ without (a).
\end{proof}

\begin{remark}\label{d:rem:volterra}
In Theorem~\ref{d:thm:volterra} the three convergences differ only through the behaviour at $s=1$, where the weight $1-s$ of $d(\cdot,\cdot)(1)$ vanishes. By (b), sequential order convergence is local $L^1$-convergence on $[0,1)$ together with a bound, a convergence of ``inductive limit'' type. It is not induced by a norm on $C[0,1]$: in a normed space every null sequence $(g_n)$ admits scalars $c_n\to\infty$ with $c_ng_n\to 0$, whereas $g_n=n^2\varphi_n\to 0$ in the sequential sense and $\nu(c_ng_n)=c_n/6\to\infty$ for every $c_n\to\infty$. This is the phenomenon behind the $C^1([0,1))$-valued fixed point examples of \cite[Theorem~3.17 and Example~6.1]{VC2}.
\end{remark}

\begin{remark}\label{d:rem:squeeze}
If $\rho(x,y)u\le d(x,y)\le C\rho(x,y)u$ for a scalar metric $\rho$ and some $u>0$, then all order convergences of $d$ reduce to $\rho$-convergence. Theorem~\ref{d:thm:band} shows that kernels positive on the diagonal behave like this even when they are not squeezed, and Theorem~\ref{d:thm:volterra} shows that kernels vanishing on the diagonal carry genuinely new information. The companion paper \cite{VC2} studies the analogous phenomenon for contractive mappings.
\end{remark}

\section{Open problems}\label{sec:open}

\begin{problem}\label{d:prob:CD}
Are \CD\ and \CDf\ equivalent for every vector lattice cover? By Proposition~\ref{d:prop:CD} this holds under positive inner density.
\end{problem}

\begin{problem}
Is \CDs\ necessary for the transfer of sequential convergence? Theorem~\ref{d:thm:seq}(iv) gives a necessary and sufficient condition in terms of sequences in $E_+$. The preorder in the proof of Theorem~\ref{d:thm:charact} produces an uncountable net even when $(b_n)$ is a sequence, so that proof does not yield a sequence $a_n\downarrow 0$. A natural intermediate condition asks that for every $b_n\downarrow 0$ in $Y$ there be a sequence $a_k\downarrow 0$ in $E$ such that every $k$ admits some $n$ with $e\le a_k$ for all $e\in E_+$ satisfying $i(e)\le b_n$; it separates the role of countability from that of positive inner density.
\end{problem}

\begin{problem}
Describe $d_k$-convergence for continuous kernels that vanish on part of the diagonal, in terms of the zero set of $s\mapsto k(s,s)$ and the order of vanishing.
\end{problem}

\begin{problem}
Define unbounded order convergence for pre-Riesz-valued vector metrics, using \cite[Definition~3.7.15 and Theorem~3.7.17]{KvG19}, and prove an analogue of Theorem~\ref{a:thm:o3}.
\end{problem}


\begin{thebibliography}{99}
\bibitem{AS05} Y.\,A. Abramovich, G. Sirotkin, On order convergence of nets, \emph{Positivity} \textbf{9} (2005), 287--292.
\bibitem{AC11} I. Altun, C. \c{C}evik, Some common fixed point theorems in vector metric spaces, \emph{Filomat} \textbf{25} (2011), 105--113.
\bibitem{Bil23} E. Bilokopytov, Order continuity and regularity on vector lattices and on lattices of continuous functions, \emph{Positivity} \textbf{27} (2023), no.~4, Paper No.~52, 29~pp.; arXiv:2103.08776.
\bibitem{Cev14} C. \c{C}evik, On continuity of functions between vector metric spaces, \emph{J. Funct. Spaces} \textbf{2014} (2014), Article ID 753969, 6~pp.
\bibitem{DC} C. \c{C}evik, Completions of vector metric spaces and the diagonal property of order convergence, preprint, arXiv:2611.12365 (2026).
\bibitem{VC2} C. \c{C}evik, Vector metrics with values in pre-Riesz spaces II: operator and nonlinear contractions, and the limits of scalarisation, preprint.
\bibitem{CA09} C. \c{C}evik, I. Altun, Vector metric spaces and some properties, \emph{Topol. Methods Nonlinear Anal.} \textbf{34} (2009), 375--382.
\bibitem{KvG19} A. Kalauch, O. van Gaans, \emph{Pre-Riesz Spaces}, De Gruyter Expositions in Mathematics 66, De Gruyter, Berlin, 2019.
\bibitem{LZ71} W.\,A.\,J. Luxemburg, A.\,C. Zaanen, \emph{Riesz Spaces I}, North-Holland, Amsterdam, 1971.
\bibitem{Roc70} R.\,T. Rockafellar, \emph{Convex Analysis}, Princeton University Press, Princeton, 1970.
\bibitem{vH93} M. van Haandel, \emph{Completions in Riesz Space Theory}, PhD thesis, University of Nijmegen, 1993.
\bibitem{Wul17} B.\,Z. Wulich, \emph{Geometrie der Kegel in normierten R\"aumen}, De Gruyter, Berlin, 2017.
\bibitem{Zab97} P.\,P. Zabrejko, $K$-metric and $K$-normed linear spaces: survey, \emph{Collect. Math.} \textbf{48} (1997), 825--859.
\end{thebibliography}
\end{document}